\documentclass[10pt]{article}
\usepackage{amssymb}
\usepackage[fleqn]{amsmath}
\usepackage{mathrsfs}
\usepackage{amsfonts}
\usepackage{tikz}
\usepackage[numbers,sort&compress]{natbib}

\usepackage{amsmath,amsthm,amsfonts,amssymb,amscd}

\allowdisplaybreaks[4]
\newtheorem {Problem} {Problem}[section]
\newtheorem {Theorem} [Problem]{Theorem}
\newtheorem {Lemma}[Problem]{Lemma}
\newtheorem{Conjecture}[Problem]{Conjecture}
\newenvironment {Proof}{\noindent {\bf Proof.}}{\hfill\ensuremath{\square}}

\begin{document}

\title{ The signless Laplacian spectral radius of graphs  without trees\thanks{This work is supported by  the National Natural Science Foundation of China (Nos.11971311, 12101166,  12161141003),  Hainan Provincial Natural Science Foundation of China (Nos. 120RC453),  Science and Technology Commission of Shanghai Municipality (No. 22JC1403602), and  the specific research fund of The Innovation Platform for Academicians of Hainan Province.
\newline \indent $^{\ddag}$E-Mails address:
mzchen@hainanu.edu.cn (M.-Z. Chen), zhaomingli@uchicago.edu  (Z.-M. Li)\ \ \ %amliu@hainanu.edu.cn (A-M. Liu),\ \ \
 xiaodong@sjtu.edu.cn (X.-D. Zhang) %(Corresponding author: X.-D.Zhang)
 }}

\author{Ming-Zhu Chen, %\thanks{E-mail: mzchen@hainanu.edu.cn},
 %A-Ming Liu,%\thanks{E-mail: amliu@hainanu.edu.cn}
 \\
School of Science, Hainan University,\\
 Haikou 570228, P. R. China, \\
\and
Zhao-Ming Li, \\
Department of Mathematics,
 The University of Chicago, \\
 %Edward H. Levi Hall
%5801 S. Ellis Ave.
Chicago, IL 60637, USA\\
%CHICAGO, IL 60615-4091
\and
 Xiao-Dong Zhang%$^{\dagger}$
\thanks{Corresponding author:  xiaodong@sjtu.edu.cn (X.-D. Zhang)
}
\\School of Mathematical Sciences, MOE-LSC, SHL-MAC\\
 % \vspace{1cm}
Shanghai Jiao Tong University,
Shanghai 200240, P. R. China
%Dedicated to Professor Fan Chung, with admiration and thanks
}

\date{}
\maketitle
\begin{abstract}
  Let $Q(G)=D(G)+A(G)$ be the signless Laplacian matrix of a simple graph of order $n$, where  $D(G)$ and $A(G)$ are the degree diagonal matrix  and the adjacency matrix  of $G$,  respectively.  In this paper, we present a sharp upper bound for the signless spectral radius of $G$ without any tree and characterize all extremal graphs which attain the upper bound, which may be regarded as a spectral extremal version for the famous Erd\H{o}s-S\'{o}s conjecture.
\\ \\
{\it AMS Subjection Classification(2020):} 05C50, 05C35\\ \\
{\it Key words:}  Signless Laplacian spectral radius;  Erd\H{o}s-S\'{o}s conjecture;    extremal graphs; tree.
\end{abstract}

\section{Introduction}
 Through this paper, we always assume that  $G$ is an undirected simple graph with vertex set
$V(G)=\{v_1,\dots,v_n\}$ and edge set $E(G)$, where $n$ is called the order of $G$.
The \emph{adjacency matrix}
$A(G)$ of $G$  is the $n\times n$  matrix $(a_{ij})$, where
$a_{ij}=1$ if $v_i$ is adjacent to $v_j$, and $0$ otherwise. The degree diagonal matrix $D(G)=diag(d_v, v\in V(G))$ of $G$ is the diagonal matrix whose diagonal elements are the degree $d_v$ of vertex $v$.   Then the  \emph{signless Laplacian spectral radius} of $G$,  denoted by $q(G)$ (or for short $q$), is the largest eigenvalue of $Q(G):=D(G)+A(G)$. Since the signless Laplacian matrix is positive definite and nonnegative, the signless Laplacian spectrum may  perform, in some sense, better  in comparison to
spectra of the adjacency matrix and other matrices associated with a graph for study of graph structure and property,    and has been extensively investigated  over the past  twenty years. For example, Cvetkovi\'{c}  and Simi\'{c} (see \cite{cvetkovic2009, cvetkovic2010-1, cvetkovic2010-2}) systematically explained a spectral theory of graphs based on the signless Laplacian matrix.      Belardo, Brunetti, Trevisan and Wang \cite{belardo2022} investigated the graphs with the signless Laplacian spectral radius at most 4.5.   Zhao, Xue and Luo \cite{zhao2020} obtained the sharp upper bound for the signless Laplacian spectral radius of $G$ in terms of the edge number and characterized all extremal graphs which attain the upper bound. For the upper and lower bounds for the signless Laplacian spectral radius are referred to \cite{DNP}, \cite{li2011}, \cite{luo2021} and references therein.    Helmberg,  Rocha and  Schwerdtfeger \cite{helmberg2017}  provided strongly polynomial time combinatorial algorithms to minimize the largest eigenvalue of  the weighted signless Laplacian of  a graph by redistributing weights among the edges.
From the view of spectral extremal graph theory,   Nikiforov \cite{Nikiforov2011} in 2011 proposed the following general problem.
\begin{Problem}
What is the maximum signless Laplacian spectral radius of a graph $G$ of
order $n$ with property $P$?
\end{Problem}
 The problem  is, in some sense,  analogue to  the core and key problem ``What is the maximum number of edges in a graph $G$ of order $n$ with
property $P$?" in extremal graph theory.  Note that this problem may be also restated to determine the maximum   largest eigenvalues  of  the matrices (such as the adjacency matrix, the signless Laplacian matrix, etc)  associated with a graph in classes of graphs,  which is known as a Brualdi-Solheid problem (see \cite{brualdi1986} ).   Due to its importance and attraction, the  problem has been developing a subject which is called spectral  extremal graph theory.  For example,  He, Jin and Zhang \cite{he2013} gave the maximum  signless Laplacian spectral radius of a graph of order $n$ without a clique of size at least $k$. Chen, Wang, Zhai \cite{chen2022} obtained the maximum signless Laplacian spectral radius of a graph of order $n$ without short cycles or long cycles.  Zhao, Huang and Guo \cite{zhao2021} determined the maximum signless Laplacian spectral radius of a graph of order $n$ forbidding a friendship graph.  Chen, Liu and Zhang \cite{Chen2020} characterized all extremal graphs having the maximum signless Laplacian spectral radius of graphs forbidding a linear forest.

 In extremal graph theory, there is a famous Erd\H{o}s-S\'{o}s conjecture which is stated as follows.
 \begin{Conjecture}\cite{erdos1965} \label{erdoscon}Let  $k$ be a positive integer and $G$ be a graph of order $n$. If the number $e(G)$ of edges  of $G$ is more than $\frac{1}{2}(k-2)n$, then $G$ contains all trees on $k$ vertices.
  \end{Conjecture}
Although there are a lot of partial results toward Conjecture~\ref{erdoscon}  (for example, see \cite{besomi2021} and references therein), it has still not been fully solved.  In spectral extremal graph theory, in 2010, Nikiforov \cite{Nikiforov2010} proposed a spectral version of the Erd\H{o}s-S\'{o}s conjecture.
 For two integers $n>k>0$, denote by $S_{n,k}$ the graph obtained by joining a clique $K_k$ with an independent set $\overline{K}_{n-k}$, and denote by $S_{n,k}^+$ the graph obtained from $S_{n,k}$ by adding one edge.

 \begin{Conjecture}\cite{Nikiforov2010}\label{nikicon} Let $k\ge 2$ and $G$ be a graph of order sufficiently large $n$ with the spectral radius $\rho(G)$ of the adjacency matrix of $G$.

 (a). If $\rho(G)\ge \rho(S_{n,k})$, then $G$ contains all trees of order $2k+2$ unless $G=S_{n,k}$;

 (b). If $\rho(G)\ge \rho(S_{n,k}^+)$, then $G$ contains all trees of order $2k+3$ unless $G=S_{n,k}^+$.
 \end{Conjecture}

In 2021, Hou, Liu,  Wang, Gao and Lv \cite{hou2021} proved that Conjecture~\ref{nikicon} is true for all trees of diameter at most 4. Very recently, Cioaba, Desai and Tait \cite{cioaba2022} confirmed conjecture~\ref{nikicon}.  It is natural to ask whether Conjecture \ref{nikicon} still holds for the signless Laplacian spectral radius of graphs. In fact,  Nikiforov and Yuan \cite{nikiforov2014} proved the following results.

 \begin{Theorem}\cite{nikiforov2014}\label{nikiforov-yuan}
 Let $k$ and $n$  be two positive integers  with $n\ge 7k^2$, and $G$ be a graph of order $n$.

(a).  If $q (G) \ge  q (S_{n,k}),$  then $G$ contains a path of order $2k+2$ unless $G = S_{n,k}$.

(b).  If $q (G) \ge  q(S^+_{n,k})$, then $G$ contains a path of order $2k+3$ unless $G = S_{n,k}^+$.
\end{Theorem}

 In this paper, motivated by their results and an approach technique from \cite{cioaba2022}, we prove the following results which extend Theorem \ref{nikiforov-yuan}.

 \begin{Theorem}\label{main}
 Let $k$  and $n$ be two positive integers with $n\ge 256000 k^8$,  and $G$ be a graph of order  $n$.

(a).  If $q (G) \ge  q (S_{n,k}),$  then $G$ contains  all trees of order $2k+2$ unless $G = S_{n,k}$.

(b).  If $q (G) \ge  q(S^+_{n,k})$, then $G$ contains  all trees of order $2k+3$ unless $G = S_{n,k}^+$.
 \end{Theorem}

  The rest of the paper is organized as follows. In Section 2, some preliminary notations and technical lemmas are presented. In Section 3, we prove Theorem \ref{main} by  progressively refining the structure of the extremal graphs having the maximum signless Laplacian spectral radius in  the set of all graphs of order $n$ not containing at least a tree on $2k+2$ vertices and  not containing at least a tree on $2k+3$ vertices, respectively.

\section{Preliminaries and lemmas}

Let $k$ be a positive integer. Denote by $\mathcal{T}_{k}$ the set of all trees on $2k+2$ vertices and $\mathcal{T}_{k}^{\prime}$ the set of all trees on $2k+3$ vertices, respectively. Moreover, denote by $\mathcal{G}_{n,k}$ the set of all graphs $G$ of order $n$ such that there exists at least one tree $T\in   \mathcal{T}_{k}$ on $2k+2$ vertices which is not a subgraph of $G$,  and $\mathcal{G}_{n,k}^{\prime}$ the set of all graphs $G$ of order $n$ such that there exists at least one tree $T\in   \mathcal{T}_{k}^{\prime}$ on $2k+3$ vertices which is not a subgraph of $G$, respectively. The rest notations and symbols may be referred to \cite{BM}.

\begin{Lemma}\label{qSnk}
Let  $k\geq2$ and $n\ge k+2$ be two integers. Then
\begin{equation}\label{qSnk-1}
n+2k-2-\frac{2k^2}{n}<q(S_{n,k})<q(S_{n,k}^+)<n+2k-2.
\end{equation}

\end{Lemma}
\begin{Proof}
From the eigenvalue and eigenvector relations, it is easy to see (for example, see \cite{nikiforov2014})
\begin{eqnarray*}
q(S_{n,k})&=&\frac{n+2k-2+\sqrt{(n+2k-2)^2-8(k^2-k)}}{2}\\
&=&n+2k-2-\frac{4(k^2-k)}{n+2k-2+\sqrt{(n+2k-2)^2-8(k^2-k)}}\\
&>& n+2k-2-\frac{4(k^2-k)}{2n}\\
&>&n+2k-2-\frac{2k^2}{n}.
\end{eqnarray*}
On the other hand, if $n=k+2$, then $q(S_{n,k}^+)=2k+2<n+2k-2$. If $n>k+2$, then similarly, by the eigenvalue and eigenvector relations, $q(S_{n,k}^+)$ is the largest root of the equation
$$z^3-(n+3k)z^2+ [(k+2)n+(2k^2-4)]z+2k^2=0.
$$
Hence $q(S_{n,k}^+)< n+2k-2$. So the assertion holds.
\end{Proof}

The following results are known fact (for example, see \cite{cioaba2022}).
\begin{Lemma}\cite{cioaba2022}\label{edgeupper} (a). If $T$ is a tree on $2k+2$ vertices and  $G$ is a graph on $n$ vertices not containing $T$ as a subgraph, then the number of edges in $G$ is at most $2kn$, i.e., $e(G)\le 2kn$.

  (b). If $T$ is a tree on $2k+3$ vertices and  $G$ is a graph on $n$ vertices not containing $T$ as a subgraph, then the number of edges in $G$ is at most $(2k+1)n$, i.e., $e(G)\le (2k+1)n$.
  \end{Lemma}

\begin{Lemma}\cite{cioaba2022} \label{bipartitecontainT}Let $t\ge 2$ be a positive integer. Then the complete bipartite graph $K_{\lfloor\frac{t}{2}\rfloor, t-1
}$ on $\lfloor\frac{t}{2}\rfloor+t-1$ vertices   contains all trees on $t$ vertices.
\end{Lemma}

In addition, for two graphs $G$ and $H$, denote by $G \vee H$ the graph obtained from two graphs $G$ and $H$ by joining all edges from each vertex in $G$ and each vertex in $H$, and denote by $\overline{G}$ the complement graph of $G$.

\begin{Lemma}\cite{cioaba2022}\label{containT2} For a positive integer $k$ and a path $P_3$ on 3 vertices,  the graph $K_{k, 2k+1}^+:= \overline{K}_k \vee ((2k-1)K_1\cup K_2)$ contains all trees on $2k+2$ vertices in $\mathcal{T}_k$, the graphs $K_{k,2k+2}^p:=\overline{K}_k \vee ((2k-1)K_1 \cup P_3)$ and $K_{k, 2k+2}^m:=\overline{K}_k \vee ((2k-2)K_1\cup 2K_2)$ contain all trees on $2k+3$ vertices in $\mathcal{T}_k^{\prime}$.
\end{Lemma}
\section{Proof of Theorem \ref{main}}

Since it is easy to see that Theorem \ref{main} holds for $k=1$, we only consider for $k\ge 2$.  Moreover, Theorem \ref{main} will hold if we can show the following theorem.

\begin{Theorem}\label{main-1}
For two positive integers $k\ge 2$ and  $n\ge 243360k^{8}$.

(a). If $G_{n,k}$ has the maximum signless Laplacian spectral radius in the set $\mathcal{G}_{n,k}$, i.e., $q(G)\le q(G_{n,k}) $ for each graph $G\in \mathcal{G}_{n,k}$,  then $G_{n,k}$ must be $S_{n,k}$.

(b). If $G_{n,k}^{\prime}$ has the maximum signless Laplacian spectral radius in the set $\mathcal{G}_{n,k}^{\prime}$, i.e., $q(G)\le q(G_{n,k}^{\prime}) $ for each graph $G\in \mathcal{G}_{n,k}^{\prime}$,  then $G_{n,k}^{\prime}$ must be $S^+_{n,k}$.
\end{Theorem}

An outline of the proof of Theorem \ref{main-1} is as follows.  The first step is to prove that both $G_{n,k}$ and $G_{n,k}^{\prime}$ are connected. The second step is to prove that there exists a vertex set of size $k$ whose vertices have  degree  at least $(1-\frac{1}{2k})n$ in $G_{n,k}$ and $G_{n,k}^{\prime}$. The third step is prove that $G_{n,k}=S_{n,k} $  and $G_{n,k}^{\prime}=S_{n,k}^+$.
% Since  proofs of  the following several Lemmas for both  $G_{n,k}=S_{n,k} $  and $G_{n,k}^{\prime}=S_{n,k}^+$ are quite similar, for brevity, we only  write all proofs for $G_{n,k}^{\prime}$.
We are now in position to prove the following lemmas.

\begin{Lemma}\label{connected}
If $n\ge \max\{2k^2, 2k+7\}$, then
 $G_{n,k}$ and $G_{n,k}^{\prime}$  are connected.
\end{Lemma}
\begin{Proof}
Suppose that $G_{n,k}^{\prime}$ is disconnected. Without loss of generality,  there exists a connected component $C_1$ such that $q(G_{n,k}^{\prime})=q(C_1)$ and there is a nonnegative eigenvector $x=(x_v, v\in V(G_{n,k}^{\prime}))$ of $G_{n,k}^{\prime}$ with maximum component 1 corresponding to the eigenvalue $q(G_{n,k}^{\prime})$ such that $Q(G_{n,k}^{\prime})x=q(G_{n,k}^{\prime})x$ with $x_z=1$ and $z\in V(C_1)$.
Let $v\in V(G_{n,k}^{\prime})\setminus V(C_1))$ and $G$ be the graph obtained from $G_{n,k}^{\prime}$ by deleting all edges incident to $v$ and adding the edge $vz$.

We have the following claim: $G\in \mathcal{G}_{n,k}^{\prime}$. In fact, note that $G_{n,k}^{\prime}\in \mathcal{G}_{n,k}^{\prime}$. There exists a tree $T\in \mathcal{T}_k$ on $2k+3$ vertices such that $G_{n,k}^{\prime}$ does not contain $T$. %So $C_1$ does not contain $T$ as a subgraph.
  Suppose that $G$ contains $T$ as a subgraph. Then $C_1+zv$ contains $T$ as a subgraph and $C_1$ does not contain $T$ as a subgraph. So $zv$ is a pendent edge in $T$ on $2k+3$ vertices.  On the other hand,  by the eigenvalue-eigenvector equation for vertex $z$,
$$q(G_{n,k}^{\prime})=  q(G_{n,k}^{\prime})x_z=d_zx_z+\sum_{u\sim z}x_u\le d_z+d_z=2d_z,$$
where $d_z$ is the degree of vertex $z$ in $G_{n,k}^{\prime}$ and  $u\sim z$ means that $u$ is adjacent to $z$.
Hence by Lemma \ref{qSnk}, $d_z\ge \frac{q(G_{n,k}^{\prime})}{2}\ge \frac{n+2k-3}{2}\ge 2k+2$, which implies that $C_1$ must contain $T$ as a subgraph.
It is a contradiction. Hence $G$ does not $T$ as a subgraph and $G\in \mathcal{G}_{n,k}^{\prime}$.
Furthermore, $q(G)\geq(C_1+zv)>q(C_1)=q(G_{n,k}^{\prime})$ which contradicts the fact that $G_{n,k}^{\prime}$ has the maximum signless Laplacian spectral radius in $\mathcal{G}_{n,k}^{\prime}$. So $G_{n,k}^{\prime}$ is connected.

The proof of the result that $G_{n,k}$ is connected is  similar and omitted.
\end{Proof}

In order to analyze the structure of $G_{n,k}$ and $G_{n,k}^{\prime}$, we introduce the following notations and symbols.
Let $\alpha$ and $\beta$ be two constants which only depend on $k$, independent of $n$, which  are  given in the following.
\begin{equation}\label{alphavalue}
\alpha=\frac{1}{80k^3}, \ \ \mbox{ and}\ \  \beta=2k\alpha.
\end{equation}
%Since  arguments for $G_{n,k}$ and $G_{n,k}^{\prime}$ are similar, we only consider the structure of $G_{n,k}^{\prime}$ except for the different place.
For the rest of paper, we always assume that $x=(x_v, v\in V(G_{n,k}^{\prime}))$ is the Perron-Frobenius vector of $Q(G_{n,k}^{\prime})$ corresponding to $q(G_{n,k}^{\prime})$ with
$$x_z=\max\{x_v, v\in V(G_{n,k}^{\prime}) \}=1.$$
    Denote by  $L$ the vertex set of $V(G_{n,k}^{\prime})$ having ``large" component value of the vector $x$ and $S$   the vertex set of $V(G_{n,k}^{\prime})$ having ``small" component value of $x$, respectively, i.e.,
\begin{equation}\label{forL}
L:=\{v\in V(G_{n,k}^{\prime})\ \ | \ x_v\ge \alpha\}, \ \ S:=\{v\in V(G_{n,k}^{\prime})\ \ | \ x_v< \alpha\}.\end{equation}
Moreover, denote by  $L^{\prime}$ the vertex set of $V(G_{n,k}^{\prime})$ having ``larger" component value of $x$ and $S^{\prime}$   the vertex set of $V(G_{n,k}^{\prime})$ having ``relative small" component value of $x$, respectively, i.e.,
\begin{equation}\label{for Lprime} L^{\prime}:=\{v\in V(G_{n,k}^{\prime})\ \ | \ x_v\ge \beta\}, \ S^{\prime}:=\{v\in V(G_{n,k}^{\prime})\ \ | \ x_v< \beta\}.
\end{equation}
In addition, let $N_i(v)$ be the vertices at distance $i$ from vertex $v$ and $L_i(v)=N_i(v)\bigcap L, $
$S_i(v)=N_i(v)\bigcap S, $   $L_i^{\prime}(v)=N_i(v)\bigcap L^{\prime}, $ and $S_i^{\prime}(v)=N_i(v)\bigcap S^{\prime} $
 for $i=1, 2$. If the vertex $v$ is unambiguous from context, we will use $L_i, S_i, L_i^{\prime}$ and $S_i^{\prime}$ instead. Moreover, $u\sim v$ means that $u$ is adjacent to $v$.  Finally, for two subsets $V_1$ and $V_2$ of $V(G_{n,k}^{\prime})$, denote by $e(V_1, V_2) $ the number of edges from  $V_1$ to $V_2$ and $|V_1|$ the size of a set $V_1$, respectively.
   Since  proofs of  the following several lemmas for both  $G_{n,k}^{\prime}=S_{n,k}^{\prime} $  and $G_{n,k}=S_{n,k}$ are quite similar, for brevity, we only  write all proofs for $G_{n,k}^{\prime}$ and the proofs for $G_{n,k}$ are omitted.
\begin{Lemma}\label{eigenvector}
Let $q$ be the signless Laplacian spectral radius of $G_{n,k}^{\prime}$ and the Perron-Frobenius vector $x=(x_v, v\in V(G_{n,k}^{\prime}))$ with $x_z=\max\{x_v, v\in V(G_{n,k}^{\prime})\}=1$. Then

\begin{equation}\label{eigenvector-1}
q^2x_v=d_v^2x_v+d_v\sum_{u\sim v}x_u+\sum_{u\sim v}d_ux_u+\sum_{u\sim v}\sum_{w\sim u}x_w,
\end{equation}

\begin{equation}\label{eigenvector-2}
d_v\sum_{u\sim v}x_u\le |L|d_v+d_v^2\alpha,
\end{equation}

\begin{equation}\label{eigenvector-3}
\sum_{u\sim v}d_ux_u\le 5kn,
\end{equation}

\begin{equation}\label{eigenvector-4}
\sum_{u\sim v}\sum_{w\sim u}x_w\le 5kn.
\end{equation}
\end{Lemma}

\begin{Proof}
It follows from the eigenvalue-eigenvector equation for $v\in V(G_{n,k}^{\prime})$ that
\begin{eqnarray*}
q^2x_v&=& q(d_vx_v+\sum_{u\sim v}x_u)\\
&=& d_v(d_vx_v+\sum_{u\sim v} x_u)+\sum_{u\sim v}(d_ux_u+\sum_{w\sim u}x_w)\\
&=& d_v^2x_v+d_v\sum_{u\sim v}x_u+\sum_{u\sim v} d_ux_u+\sum_{u\sim v}\sum_{w\sim u}x_w.
\end{eqnarray*}
So (\ref{eigenvector-1}) holds. By the definition of $L$ and $S$,
$$\sum_{u\sim v}x_u=\sum_{u\sim v, u\in L}x_u+\sum_{u\sim v, u\in S}x_u\le \sum_{u\sim v, u\in L}1+\sum_{u\sim v, u\in S}\alpha\le |L|+d_v\alpha.$$
So (\ref{eigenvector-2}) holds.
By Lemma \ref{edgeupper} (b), we have $e(G_{n,k}^{\prime})\le (2k+1)n$ and
\begin{eqnarray*}
\sum_{u\sim v}d_ux_u\le \sum_{u\sim v}d_u\le  2e(G_{n,k}^{\prime})\leq(4k+2)n\leq 5kn.
\end{eqnarray*}
So (\ref{eigenvector-3}) holds. Moreover,

$$\sum_{u\sim v}\sum_{w\sim u}x_w\le \sum_{u\sim v}d_u\le2e(G_{n,k}^{\prime})\leq(4k+2)n\leq 5kn.$$
So (\ref{eigenvector-4}) holds.
\end{Proof}

\begin{Lemma}\label{Lupper} Let $L:=\{v\in V(G_{n,k}^{\prime})\  | \ x_v\ge \alpha\}.$  Then
$|L|\le \frac{10k}{\alpha}$.
\end{Lemma}

\begin{Proof}   By Lemma \ref{edgeupper} (b),  $\sum_{v\in L}d_v\le 2e(G_{n,k}^{\prime})\le (4k+2)n\leq 5kn$.
Then by the definition of $L$,
$$q|L|\alpha\le \sum_{v\in L}qx_v=\sum_{v\in L}\bigg(d_vx_v+\sum_{u\sim v}x_u\bigg)\le \sum_{v\in L}2d_v\le 10kn.  $$
Hence by  (\ref{qSnk-1}) in  Lemma \ref{qSnk},
\begin{equation}\label{Lupper-2}
|L|\le \frac{10kn}{q\alpha}\le \frac{10k}{\alpha}.\nonumber
\end{equation}
\end{Proof}
%\begin{Lemma}\label{Lupper}
%If $n\ge \frac{64k}{\alpha^3}$  then
%\begin{equation}\label{Lupper-1}
%|L|\le 7k.
%\end{equation}
%\end{Lemma}
%
%\begin{Proof}  Notice that $L=\{v\in V(G_{n,k}\ \ |\ x_v\ge \alpha\}$ and $\sum_{v\in L}d_v\le 2e(G_{n,k})\le 4kn$ by Lemma \ref{edgeupper}.
%Then
%$$q|L|\alpha\le \sum_{v\in L}qx_v=\sum_{v\in L}(d_vx_v+\sum_{u\sim v}x_u)\le \sum_{v\in L}2d_v\le 8kn.  $$
%Hence by Lemma \ref{qSnk},
%\begin{equation}\label{Lupper-2}
%|L|\le \frac{8kn}{q\alpha}\le \frac{8k}{\alpha}.
%\end{equation}
%Furthermore, by $(q-d_v)x_v=\sum_{u\sim v}x_u$, $e(L,s)\le e(G_{n,k})\le 2kn$,  (\ref{Lupper-2}), and $\left(\frac{8k}{\alpha}\right)^2\le kn\alpha$ from $n\ge  \frac{64k}{\alpha^3}$,
%\begin{eqnarray*}
%\sum_{v\in L}(q-d_v)\alpha &\le & \sum_{v\in L}(q-d_v)x_v\\
%&=& \sum_{v\in L}\sum_{u\sim v, u\in S}x_u+\sum_{v\in L}\sum_{u\sim v, u\in L}x_u\\
%&\le&  e(L, S) \alpha+|L|^2\\
%&\le & 4kn\alpha+\left(\frac{8k}{\alpha}\right)^2\\
%&\le& 5kn \alpha.
%\end{eqnarray*}
%Hence by $\sum_{v\in L}d_v\le 2e(G_{n,k})\le 4kn$ and
%$$q|L|\alpha\le \sum_{v\in L}d_v\alpha+ 5kn\alpha\le 7kn\alpha.$$
%So $$|L|\le \frac{7kn\alpha}{q\alpha}\le 7k.$$
%\end{Proof}
\begin{Lemma}\label{degree}
Let $L^{\prime}=\{v\in V(G_{n,k}) \  | \ x_v\ge\beta \}$  with $\beta=2k\alpha$. If $n\geq\frac{40k^2}{\alpha^2}$, then
$d_v\ge (1-\frac{1}{2k})n$ for all $v\in L^{\prime}$.
\end{Lemma}
\begin{Proof}
Suppose that there would exist a vertex $s\in L^{\prime}$ such that $d(s)<(1-\frac{1}{2k})n$. Then
by (\ref{eigenvector-1})--(\ref{eigenvector-4}) in  Lemma \ref{eigenvector}  and $|L|\le \frac{10k}{\alpha}$ in  Lemma \ref{Lupper},
\begin{eqnarray*}
q^2x_s&=&d_s^2x_+d_s\sum_{u\sim s}x_u+\sum_{u\sim s}d_ux_u+\sum_{u\sim s}\sum_{w\sim u}x_w\\
&\le & d_s^2(x_s+\alpha)+d_s|L|+5kn+5kn\\
%using $|L|\le 7k$  (need $n>\frac{64k}{\alpha^3}$
&\le & d_s^2(x_s+\alpha)+\frac{10kn}{\alpha}+10kn\\
&= & d_s^2(x_s+\alpha)+\bigg(10+\frac{10}{\alpha}\bigg)kn\\
&\le & (x_s+\alpha)\bigg(1-\frac{1}{2k}\bigg)n^2+\bigg(10+\frac{10}{\alpha}\bigg)kn.
\end{eqnarray*}
Hence by $x_s\ge \beta=2k\alpha$, $n\ge \frac{40k^2}{\alpha^2}$ and Lemma \ref{qSnk},
\begin{eqnarray*}
q^2&\le & \bigg(1+\frac{\alpha}{x_s}\bigg)\bigg(1-\frac{1}{2k}\bigg)n^2+\bigg(10+\frac{10}{\alpha}\bigg)\frac{kn}{x_s}\\
 &\le & \bigg(1+\frac{1}{2k}\bigg)\bigg(1-\frac{1}{2k}\bigg)n^2+\frac{10n}{\alpha^2}\\
 &= & n^2-\left(\frac{n}{4k^2}-\frac{10}{\alpha^2}\right)n\\
 %here $n\ge frac{19k}{\alpha}
 &\le & n^2 < q^2,
 \end{eqnarray*}
 which is a contradiction. Hence the assertion holds.
\end{Proof}
 \begin{Lemma}\label{Lprimelek}
Let $L^{\prime}=\{v\in V(G_{n,k}) \  |\ x_v\ge\beta\ \}$  with $\beta=2k\alpha$. If $n\geq\frac{40k^2}{\alpha^2}$, then
$|L^{\prime}|\le k$.
\end{Lemma}
\begin{Proof}
If $|L^{\prime}|\geq k+1$, then we choose arbitrarily $k+1$ distinct vertices $v_1, v_2,\dots, v_{k+1}\in L^{\prime}$. By Lemma~\ref{degree},
$$\left|\bigcap_{i=1}^{k+1} N(v_i)\right|\geq \sum_{i=1}^{k+1} d_{v_i}-k\bigg|\bigcup_{i=1}^{k+1} N(v_i)\bigg|=(k+1)\bigg(1-\frac{1}{2k}\bigg)n-kn %\bigg(1-\frac{k+1}{2k}\bigg)n\geq\frac{n}{4}
\geq 2k+2.$$
So $G_{n,k}^{\prime}$ contains a complete bipartite graph $ K_{k+1,2k+2}$ as a subgraph. By Lemma~\ref{bipartitecontainT}, $G_{n,k}^{\prime}$ contains all trees on $2k+3$ vertices, which is a contradiction. Hence the assertion holds.
\end{Proof}

It is ready to determine the exact value of $|L^{\prime}|$.

\begin{Lemma}\label{Lprimekk} If $n\ge \frac{40k^2}{\alpha^2}$, %and $\alpha\le \frac{1}{26k^2}$,
then $|L^{\prime}|=k$.
\end{Lemma}

\begin{Proof}
 By Lemma~\ref{Lprimelek}, $|L^{\prime}|\leq k$. Assume for a contradiction that $|L^{\prime}|\le k-1$. By the eigenvalue-eigenvector equation for the vertex $z$, we have
 \begin{equation}\label{Lprimekk-0}
 q(q-d_z)=q(q-d_z)x_z=\sum_{u\sim z}qx_u=\sum_{u\sim z}d_ux_u+\sum_{u\sim z}\sum_{w\sim u}x_w.
 \end{equation}
 Recalling the definition of $L_i^{\prime}=N_i(z)\bigcap L^{\prime}$ and $S_i^{\prime}=N_i(z)\bigcap S^{\prime}$ for $i=1,2$,  where
 $L^{\prime}=\{v\in V(G_{n,k}^{\prime}) \ \ | \ x_v\ge \beta\ \}$ and $S^{\prime}=\{v\in V(G_{n,k}^{\prime}) \ \ | \ x_v< \beta\ \}$  for $\beta=2k\alpha$.
 Then $ L_1^{\prime}\bigcup L_2^{\prime}\subseteq L^{\prime}\setminus \{z\}$ and
 \begin{eqnarray}\label{Lprimekk-1}
 \sum_{u\sim z}d_ux_u &=& \sum_{u\sim z, u\in S_1^{\prime}}d_ux_u+ \sum_{u\sim z, u\in L_1^{\prime}}d_ux_u\nonumber\\
 &\le & \beta  \sum_{u\sim z, u\in S_1^{\prime}}d_u +\sum_{u\sim z, u\in L_1^{\prime}}n \nonumber\\
 &\le & 2e(G_{n,k}^{\prime})\beta+ (|L^{\prime}|-1)n \nonumber\\
 &\le & 5kn\beta+(k-2)n.
 \end{eqnarray}
  Furthermore, by Lemma \ref{edgeupper} (b), we have%,  and $|L^{\prime}|\le k$ in Lemma \ref{Lprimelek},
 \begin{eqnarray}\label{Lprimekk-2}
 \sum_{u\sim z}\sum_{w\sim u}x_w %\nonumber\\
&=& \sum_{u\sim z, u\in S^{\prime}}\sum_{w\sim u}x_w+\sum_{u\sim z, u\in L^{\prime}}\sum_{w\sim u}x_w \nonumber \\
&=&\sum_{ u\in S_1^{\prime}}\left(x_z+ \sum_{w\sim u, w\in S_1^{\prime} \cup S_2^{\prime}}x_w
+\sum_{w\sim u, w\in L_1^{\prime} \cup L_2^{\prime}}x_w\right)+ \nonumber\\
&& \sum_{ u\in L_1^{\prime}}\left(x_z+ \sum_{w\sim u, w\in S_1^{\prime} \cup S_2^{\prime}}x_w
+\sum_{w\sim u, w\in L_1^{\prime} \cup L_2^{\prime}}x_w\right) \nonumber\\
& \le & d_z+e(S_1^{\prime},  S_1^{\prime} \cup S_2^{\prime})\beta + e(S_1^{\prime}, L_1^{\prime}\cup L_2^{\prime})
+e(L_1^{\prime}, S_1^{\prime}\cup S_2^{\prime})\beta + e(L_1^{\prime}, L_1^{\prime}\cup L_2^{\prime})\nonumber\\
&\le & d_z+ 2e(G_{n,k}^{\prime})\beta+|L_1^{\prime}\cup L_2^{\prime}|n+e(G_{n,k}^{\prime})\beta+|L|^2 \nonumber\\
&\le &  (k-1)n+ 8kn\beta+(k-1)^2.
\end{eqnarray}
%Hence by (\ref{eigenvector-1}) and (\label{Lprimekk-2}), and $13k\beta=26k^2\alpha<1$, and $n>7k^2>k^2+6k$,  we have
With  $13k\beta=26k^2\alpha<1$,  (\ref{Lprimekk-1}) and (\ref{Lprimekk-2}) yield
 \begin{equation}\label{Lprimekk-3}
 \sum_{u\sim z}d_ux_u+\sum_{u\sim z}\sum_{w\sim u}x_w\le (2k-3)n+13k\beta n+(k-1)^2<(2k-2)n+(k-1)^2.
 \end{equation}
By Lemma \ref{qSnk} and $n\ge \frac{40k^2}{\alpha^2}\geq8k-6$, we have
\begin{eqnarray}\label{Lprimekk-4}
(q-d_z)q&\ge & \left (n+2k-2-\frac{2k^2}{n}-(n-1)\right)\left (n+2k-2-\frac{2k^2}{n}\right)\nonumber\\
&=& (2k-1)n+(2k-1)\bigg(2k-2-\frac{2k^2}{n}\bigg)-\frac{2k^2}{n}\bigg(n+2k-2-\frac{2k^2}{n}\bigg)\nonumber\\
&=&(2k-1)n+2k^2-6k+2-\frac{2k^2(4k-3)}{n}+\frac{4k^2}{n^2}\nonumber\\
&>&(2k-1)n+k^2-6k+2+k^2-\frac{2k^2(4k-3)}{n}\nonumber\\
&\ge & (2k-1)n+k^2-6k+2
\end{eqnarray}
Therefore, substituting (\ref{Lprimekk-3}) and (\ref{Lprimekk-4}) into (\ref{Lprimekk-0}), we have
$$(2k-1)n+k^2-6k+2\le (q-d_z)q=(q-d_z)qx_z<(2k-2)n+(k-1)^2, $$
%By some algebra, we have  $n<4k-1$,
which is a contradiction. Hence $|L^{\prime}|= k$.
\end{Proof}

\begin{Lemma}\label{eigenvectorupper}
If $n\geq\frac{40k^2}{\alpha^2}$, % and $\alpha\le \frac{1}{78k^3}$,
 then $x_v\ge 1-\frac{1}{k}$ for all $v\in L^{\prime}$.
 \end{Lemma}
 \begin{Proof}
 The assertion is proved by contradiction.  Suppose that there would  exist a vertex $v\in L^{\prime}$ with $x_v<1-\frac{1}{k}$.
 By the eigenvalue-eigenvector equation for the vertex $z$,
 \begin{equation}\label{eigenvectorupper-1}
  q(q-d_z)=q(q-d_z)x_z=\sum_{u\sim z}qx_u=\sum_{u\sim z}d_ux_u+\sum_{u\sim z}\sum_{w\sim u}x_w.
 \end{equation}
 % Notice that
 \begin{eqnarray}\label{eigenvectorupper-2}
 \sum_{u\sim z}d_ux_u &=& \sum_{u\sim z, u\in S_1^{\prime}}d_ux_u+ \sum_{u\sim z, u\in L_1^{\prime}}d_ux_u\nonumber\\
 &\le & \beta  \sum_{u\sim z, u\in S_1^{\prime}}d_u +\sum_{u\sim z, u\in L_1^{\prime}}n \nonumber\\
 &\le & 2e(G_{n,k}^{\prime})\beta+ (|L^{\prime}|-1)n \nonumber\\
 &\le & 5kn\beta+(k-1)n.
 \end{eqnarray}
% Denote by $N_{S_1^{\prime}}(z)=N(z)\bigcap S_1$,
By Lemma \ref{degree}, $|N(v)\cap N(z)|=|N(v)|+|N(z)|-|N(v)\cup N(z)|\ge 2(1-\frac{1}{2k})n-n=(1-\frac{1}{k})n>0$. Hence $v\in N_1(z)\cup N_2(z)$ and $v\in L_1^{\prime}\cup L_2^{\prime}$.  %Moreover, by Lemma \ref{degree}, $|N(z)\cap N(v)|\ge |N(z)|+|N(v)|-n\ge (1-\frac{2}{k^2})n$.
So\begin{equation}\label{eigenvectorupper-3}
e(S_1^{\prime}, v)=  |N(v)\cap N(z)\cap S^{\prime}|=|N(v)\cap N(z)|-|N(v)\cap  N(z)\cap L^{\prime}|\ge \bigg(1-\frac{1}{k}\bigg)n -k.
\end{equation}
Then by (\ref{eigenvectorupper-3}), Lemmas \ref{edgeupper} (b)  and  \ref{Lprimekk}, we have
  \begin{eqnarray}\label{eigenvectorupper-4}
&& \sum_{u\sim z}\sum_{w\sim u}x_w \nonumber\\
&=& \sum_{u\sim z, u\in S^{\prime}}\sum_{w\sim u}x_w+\sum_{u\sim z, u\in L^{\prime}}\sum_{w\sim u}x_w \nonumber \\
&=&\sum_{ u\in S_1^{\prime}}\left(x_z+ \sum_{w\sim u, w\in S_1^{\prime} \cup S_2^{\prime}}x_w
+\sum_{w\sim u, w\in L_1^{\prime} \cup L_2^{\prime}\setminus\{v\}}x_w\right)+ \sum_{ u\in S_1^{\prime}}\sum_{v\sim u}x_v+\nonumber\\
&& \sum_{ u\in L_1^{\prime}}\left(x_z+ \sum_{w\sim u, w\in S_1^{\prime} \cup S_2^{\prime}}x_w
+\sum_{w\sim u, w\in L_1^{\prime} \cup L_2^{\prime}}x_w\right) \nonumber\\
& \le & d_z+e(S_1^{\prime},  S_1^{\prime} \cup S_2^{\prime})\beta + e(S_1^{\prime}, L_1^{\prime}\cup L_2^{\prime}\setminus\{v\})+
| N(v)\cap S_1^{\prime}|x_v+\nonumber\\
&&e(L_1^{\prime}, S_1^{\prime}\cup S_2^{\prime})\beta + e(L_1^{\prime}, L_1^{\prime}\cup L_2^{\prime})\nonumber\\
&\le & d_z+ 2e(G_{n,k}^{\prime})\beta+ e(S_1^{\prime}, L_1^{\prime}\cup L_2^{\prime})-e(S_1^{\prime}, v) +  e(S_1^{\prime}, v)\bigg(1-\frac{1}{k}\bigg)+e(G_{n,k}^{\prime})\beta+|L^{\prime}|^2 \nonumber\\
&\le &  kn+ 8kn\beta-\bigg(\bigg(1-\frac{1}{k}\bigg)n-k\bigg)\frac{1}{k}+k^2.\nonumber\\
&= &  kn+ 8kn\beta-\frac{(k-1)n}{k^2}+k^2+1.
\end{eqnarray}
Hence  by $n\ge \frac{40k^2}{\alpha^2}\geq 18k^3$ and  $13k\beta=26k^2\alpha\leq\frac{1}{3k}$, (\ref{eigenvectorupper-2})  and (\ref{eigenvectorupper-4}) yield
\begin{eqnarray}\label{eigenvector-5}
&& \sum_{u\sim z}d_ux_u +\sum_{u\sim z}\sum_{w\sim u}x_w \nonumber\\
&\le & (2k-1)n+13k\beta n -\frac{k-1}{k^2}n+k^2+1\nonumber\\
%&=& (2k-1)n+26k\alpha n -\frac{k-1}{k^2}n+k^2+1\nonumber\\
&\leq& (2k-1)n+\left(\frac{1}{3k}-\frac{1}{k}+\frac{1}{k^2}\right)n+k^2+1\nonumber\\
%&\le & (2k-1)n-\frac{(2k-3)n}{3k^2}+k^2 +1\nonumber\\
&\le & (2k-1)n-\frac{n}{3k}+k^2 +1\nonumber\\
&\le & (2k-1)n+k^2-6k+2.
\end{eqnarray}
%On the other hand,
%by Lemma \ref{qSnk} and,%$n\ge \frac{40k^2}{\alpha^2}\geq8k-6$, we have
%\begin{eqnarray}\label{eigenvector-6}
%(q-d_z)q&\ge & \left (n+2k-2-\frac{2k^2}{n}-(n-1)\right)\left (n+2k-2-\frac{2k^2}{n}\right)\nonumber\\
%&=& (2k-1)n+(2k-1)\bigg(2k-2-\frac{2k^2}{n}\bigg)-\frac{2k^2}{n}\bigg(n+2k-2-\frac{2k^2}{n}\bigg)\nonumber\\
%&=&(2k-1)n+2k^2-6k+2-\frac{2k^2(4k-3)}{n}+\frac{4k^2}{n^2}\nonumber\\
%&>&(2k-1)n+k^2-6k+2+k^2-\frac{2k^2(4k-3)}{n}\nonumber\\
%&\ge & (2k-1)n+k^2-6k+2.
%\end{eqnarray}
Hence by  (\ref{Lprimekk-4}) in Lemma~\ref{Lprimekk},  (\ref{eigenvectorupper-1}) and (\ref{eigenvector-5}), we have
$$(2k-1)n+k^2-6k+2<(q-d_z)q=\sum_{u\sim z}d_ux_u+\sum_{u\sim z}\sum_{w\sim u}x_w\leq (2k-1)n+k^2-6k+2, $$
 which is a contradiction. So the assertion holds.
  \end{Proof}

\begin{Lemma} \label{structure}
Let $L^{\prime}=\{v\in V(G_{n,k}^{\prime})\ |\ x_v\ge \beta\ \}$ and $R=\bigcap_{v\in L^{\prime}}N(v)$ with $\beta=2k\alpha$.
If $n\geq\frac{40k^2}{\alpha^2}$, %and $\alpha \leq\frac{1}{78k^3}$,
then $|L^{\prime}|=k$ and $|R|=n-k$.
\end{Lemma}
\begin{Proof}
It follows from Lemma \ref{Lprimekk} that $|L^{\prime}|=k$.
By Lemma \ref{degree},
$$|R|= |\bigcap_{v\in L^{\prime}}N(v)|\ge \sum_{v\in L^{\prime}}|N(v)|-(k-1)n\ge k\bigg(1-\frac{1}{2k}\bigg)n-(k-1)n= \frac{n}{2}.$$
We now prove that $|R|=n-k$ by contradiction.  Suppose that
$F=V(G_{n,k}^{\prime})\setminus (L^{\prime}\cup R)\neq \emptyset$.  Noticing  that $G_{n,k}^{\prime}\in \mathcal{G}_{n,k}^{\prime}$,  there exists a tree $T$ on  $2k+3$ vertices such that $G_{n,k}^{\prime}$ does not contain $T$ as a subgraph. Then the following (a)-(d) hold.

 (a).  For each $v\in F$, $e(v, R)\le 2k+1$. Otherwise the bipartite graph consisting of two parts $L^{\prime}\cup \{v\}$ and $R$ contains a complete bipartite graph $K_{k+1, 2k+2}$, which implies that $G_{n,k}^{\prime}$ contains all trees on $2k+3$ vertices. This is a contradiction.

 (b). $ e(F)\le (2k+1) |F|$.  Since $G_{n,k}^{\prime}$ does not contains some tree $T$ on $2k+3$ vertices which implies that the subgraph od $G_{n,k}^{\prime}$ by induced by $F$ does not contain $T$. Hence by Lemma \ref{edgeupper} (b), $ e(F)\le (2k+1) |F|$.

 (c). There exists a vertex $s\in F$ such that the degree of $s$ in $F$ at most $5k$.  The assertion follows from  (b).

(d).  For $s\in F$ in (c),  $e(s, L^{\prime})\le k-1$. In fact,  $s\notin R$  implies that $s$ is not adjacent to all vertices in $L^{\prime}$.  So
 $e(s, L^{\prime})\le |L^{\prime}|-1=k-1.$ Moreover, there exists a vertex $u\in L^{\prime}$ such that $s$ is not adjacent to $u$.

 Let $G$ be the graph obtained from $G_{n,k}^{\prime}$ by deleting all edges $sw$ incident to $s$  for $w\in R\cup F$ and  adding one edge $su$.
Then it is easy to see that $G$ does not contain $T$ as a subgraph. Furthermore, by $6k\beta=12k^2\alpha\leq \frac{12k^2}{78k^3}=\frac{2}{13k}$,
\begin{eqnarray*}\label{structure-1}
&&x^TQ(G)x-x^TQ(G_{n,k})x\nonumber\\
&=& (x_s+x_u)^2-\sum_{v\sim s, v\in R\cup F}(x_v+x_s)^2\nonumber\\
&\ge &\bigg (1-\frac{1}{k}\bigg)^2-(5k+2k+1)4\beta^2\nonumber\\
&\ge & \bigg(1-\frac{1}{k}\bigg)^2-(6k\beta)^2\nonumber\\
%&\ge & (1-\frac{1}{k}+6k\beta)(1-\frac{1}{k}-6k\beta)\nonumber\\
%&\ge & (1-\frac{1}{k}+6k\beta)\bigg(1-\frac{1}{k}-\frac{2}{13k}\bigg)\nonumber\\
%&= &(1-\frac{1}{k}+6k\beta)\bigg(1-\frac{15}{13k}\bigg)\nonumber\\
&> & 0.
\end{eqnarray*}
% Since $G_{n,k}\in \mathcal{G}_{n,k}$, there exists a tree $T$ on
So $q(G)>q(G_{n,k}^{\prime})$  and $G\in \mathcal{G}_{n,k}^{\prime}$ which contradicts that $G_{n,k}^{\prime}$ has the maximum signless Laplacian spectral radius in $\mathcal{G}_{n,k}^{\prime}$. Hence $F=\emptyset$ and $|R|=n-k$.
\end{Proof}

 Now we are ready to prove Theorem \ref{main-1}.

 {\bf Proof} (a). Let $G_{n,k}$ be the maximum signless spectral radius of graph in $\mathcal{G}_{n,k}$.  Then by Lemma \ref{structure}, the vertex set of $G_{n,k}$  consists of $L^{\prime}$ with size $k$ and $R$ with size $n-k$ such that each vertex in $L^{\prime}$ is adjacent to each vertex in $R$.  By Lemma \ref{containT2}, $R$ must be an independent set.   So $G_{n,k}$ is a subgraph of $S_{n,k}$  which implies that $q(G_{n,k})\le q(S_{n,k})$. On the other hand, $S_{n,k}\in \mathcal{G}_{n,k}$ which implies that $q(S_{n,k})\le q(G_{n,k})$.  So  $q(G_{n,k})=q(S_{n,k})$ and $G_{n,k}$ has to be $S_{n,k}$.

 (b). Let $G_{n,k}^{\prime}$ be the maximum signless spectral radius of graph in $\mathcal{G}_{n,k}^{\prime}$.  Then by Lemma \ref{structure}, the vertex set of $G_{n,k}^{\prime}$ consists of $L^{\prime}$ with size $k$ and $R$ with size $n-k$ such that each vertex in $L^{\prime}$ is adjacent to each vertex in $R$. By Lemma \ref{containT2}, $R$ contains at most one edge.    So $G_{n,k}^{\prime}$ is a subgraph of $S_{n,k}^+$  which implies that $q(G_{n,k}^{\prime})\le q(S_{n,k}^+)$. On the other hand, $S_{n,k}^+\in \mathcal{G}_{n,k}^{\prime} $ which implies that $q(S_{n,k}^+)\le q(G_{n,k}^{\prime})$.  So  $q(G_{n,k}^{\prime})=q(S_{n,k}^+)$ and $G_{n,k}^{\prime} $ has to be $S_{n,k}^+$.

%\subsection*{Acknowledgements}
%The authors  are grateful to the  referee for him/her helpful comments  and suggestion of the manuscript.

\end{document}